\documentclass[12pt]{amsart}
\usepackage[applemac]{inputenc}
\usepackage{bm}

\usepackage[dvipsnames,svgnames,x11names,hyperref]{xcolor}
\usepackage[hyphens]{url}
\usepackage[colorlinks=true,linkcolor=Maroon,citecolor=blue,urlcolor=blue,hypertexnames=false,linktocpage]{hyperref}
\usepackage{bookmark}
\usepackage{amsmath,thmtools,mathtools}
\mathtoolsset{showonlyrefs=true}

\usepackage{blindtext}

\newcounter{mgncount}

\usepackage{fancyhdr}
\usepackage{esint}
\usepackage{enumerate}

\usepackage{pictexwd,dcpic}
\usepackage{graphicx}
\usepackage{caption}
\usepackage{graphicx}

\declaretheorem[name=Theorem,numberwithin=section]{thm}

\declaretheorem[name=Lemma,sibling=thm]{lemma}

\declaretheorem[name=Corollary,sibling=thm]{cor}

\declaretheorem[style=remark,name=Remark,numbered=no]{remark}

\numberwithin{equation}{section}

\newcommand{\bbR}{\mathbb{R}}

\newcommand{\bbS}{\mathbb{S}}

\newcommand{\Sn}{\mathbb{S}^n}

\newcommand{\cH}{\mathcal{H}}

\newcommand{\cK}{\mathcal{K}}

\newcommand{\pf}[1]{\begin{proof}#1 \end{proof}}
\newcommand{\eq}[1]{\begin{equation}\begin{alignedat}{2} #1 \end{alignedat}\end{equation}}

\newcommand{\hp}{\hphantom}

\begin{document}
\title[]
 {Uniqueness of solutions to a class of non-homogeneous curvature problems}
\author[]{Mohammad N. Ivaki}
\maketitle

\begin{abstract}
We show that the only even, smooth, convex solutions to a class of isotropic mixed Christoffel-Minkowski type problems are origin-centred spheres, which, in particular, answers a question of Firey 74 in the even isotropic case about kinematic measures. Employing the Heintze-Karcher inequality, we prove that the only smooth, strictly convex solutions to a large class of Minkowski type problems are origin-centred spheres. Immediate corollaries are the uniqueness of solutions to the isotropic Orlicz-Minkowski problem and the isotropic $L_p$-Gaussian-Minkowski problem when $p\geq 1$.
\end{abstract}

\section{Introduction}

Let $(\bbR^{n+1},\delta=\langle\,,\rangle,D)$ be the Euclidean space with its standard inner product and flat connection, and
let $(\bbS^n,\bar{g},\bar{\nabla} )$ denote the unit sphere equipped with its standard round metric and Levi-Civita connection. 
 
Suppose $f:\bbS^n\to (0,\infty)$ and $\psi: (0,\infty)\times (0,\infty)\to (0,\infty)$ are smooth functions. Consider the curvature equation
\eq{\label{main eq}
\sum_k \alpha_k\sigma_k(\bar{\nabla}^2u+\bar{g}u,\bar{g})=\psi(u,|Du|)f.
}
Here $\sigma_k=\sigma_k(\bar{\nabla}^2u+\bar{g}u,\bar{g})$ is $k$th elementary symmetric function of eigenvalues of $\bar{\nabla}^2u+\bar{g}u$ with respect to $\bar{g}$, $\alpha_k$ are non-negative constants, and at least two of them are positive. The left-hand side of \eqref{main eq} appears as densities of the kinematic measures on convex bodies. When $\psi\equiv 1$, uniqueness and existence questions for kinematic measures (in the class of convex bodies) were posed by Firey \cite{Fir74} and Schneider \cite{Sch76}. See also \cite[p. 454]{Sch14}. For $n=2$, the uniqueness (up-to translations) of sufficiently smooth convex solutions to \eqref{main eq} with $\psi\equiv1$ follows from a theorem of A. D. Aleksandrov (see \cite[Sec.~9]{Fir74} and \cite{GWZ16}), and for general measures it was proved in \cite{Sch76}. The existence for $n\geq 2$ and uniqueness for $n\geq 3$ are still open. See also \cite{GZ21} and \cite[Thm.~2.4]{BIS23}.

Here, as our first result, we prove the uniqueness of even, smooth, convex solutions to \eqref{main eq} in the isotropic case (i.e. $f\equiv 1$), which, in particular, answers Firey's question in this case. The reader may consult \cite{Lut93, And99, BCD17, Che20,Sar22, LW22, IM23} and the references therein for uniqueness results when only one $\alpha_k$ is non-zero; in particular, concerning the Gauss curvature type problems.

We say $u\in C^2(\bbS^n)$ is convex if the eigenvalues of $\bar{\nabla}^2u+\bar{g}u$ with respect to $\bar{g}$ are all non-negative. Moreover, $u$ is said to be even if $u(x)=u(-x)$ for all $x\in\Sn$.

\begin{thm}\label{sigma k soliton}
Suppose $\varphi: (0,\infty)\times (0,\infty)\to (0,\infty)$ is $C^1$-smooth, $\alpha_k\geq 0$ and at least two of them are non-zero. Let  $0<u\in C^2(\bbS^n)$ be an even, convex solution of
\eq{\label{main eq1}
\sum_k \alpha_ku\sigma_k=\varphi(u,|Du|).
}
If $\ell-1+x\partial_1(\log\varphi)(x,y)> 0$ and $\partial_2\varphi\geq 0$, where $\ell:=\min \{k: \alpha_k>0\}$,
then $u$ is constant.
\end{thm}

As a complementary result to \autoref{sigma k soliton}, let us mention that Kohlmann \cite{Koh98}, using the Heintze-Karcher inequality, showed if a non-negative linear combination of curvature measures of a convex body is proportional to the boundary measure, then the convex body is a ball. Recently Kohlman's theorem was extended to the anisotropic setting in \cite{ALWX21}. Here we also use the Heintze-Karcher inequality to prove the following new uniqueness result.

\begin{thm}\label{Gauss1}
Suppose $\varphi: (0,\infty)\times(0,\infty)\to (0,\infty)$ is $C^1$-smooth with $\partial_1\varphi\geq0$, $\partial_2\varphi\geq 0$ and at least one of these inequalities is strict.
If $M^n$ is a closed, smooth, strictly convex hypersurface with the support function $h>0$ and Gauss curvature $\cK$, such that
\eq{
\varphi(h,|Dh|)\mathcal{K}=1,
}
then $M^n$ is a rescaling of $\bbS^n$. Moreover, if $\partial_1\varphi=0$, then the same conclusion holds under no extra assumption on the sign of $h$. 
\end{thm}

The theorem implies the uniqueness of solutions to the isotropic Orlicz-Minkowski problem, i.e. $\partial_2\varphi=0$, (the corresponding $L_p$ case is $p>1$), which was introduced by Chou, Wang \cite{CW06} and Haberl, Lutwak, Yang, Zhang in \cite{HLYZ10}. It is worth pointing out that while for the case $\varphi(h)=h^{p-1}$ and $p>1$, the uniqueness follows immediately from the $L_p$-Minkowski inequality \cite{Lut93}, in the non-homogenous case, the Orlicz-Minkowski inequality (cf. \cite{XJL14,GHW14,GHWXY20,GHWXY19}) by itself does not imply uniqueness. See also \cite{Sar22} for a stronger result when $\partial_2\varphi=0$.

Another corollary of \autoref{Gauss1} is a uniqueness result in connection to the $L_p$-Gaussian-Minkowski problem, which affirmatively answers a conjecture by Chen, Hu, Liu, and Zhao in \cite{CHLZ23}. 

\begin{cor}\label{lp gaussian}
Let $p>1$ and $M^n$ be a closed, smooth, strictly convex hypersurface with the support function $h$ and Gauss curvature $\cK$. If $h>0$ and
$
h^{p-1}e^{\frac{1}{2}|Dh|^2}\mathcal{K}=c,
$
then $M^n$ is a rescaling of $\bbS^n$. The same conclusion holds for $p=1$ under no extra assumption on the sign of $h$.
\end{cor}

For the $n=1, p=1, h\geq0$ case, the corollary was first proved in \cite{CHLZ23} based on an argument of Andrews \cite{And03}. Assuming that $M^n$ is origin-centred but allowing $p>-n-1$, the corollary was proved in \cite{IM23} via the local Brunn-Minkowski inequality. It is of great interest to see whether for $p<1$ the origin-centred assumption could be removed, perhaps by substituting suitable test functions into the local Brunn-Minkowski inequality.

The strength of \autoref{lp gaussian} stems from the fact that we make no additional assumption on the size of the constant $c$. Such uniqueness results are key to obtaining existence results via the degree theory, hence allowing the elimination of the Lagrange multiplier that emerges from variational and flow approaches. Under an additional assumption on the size of $c$, the uniqueness follows from the Ehrhard inequality and the characterization of its equality cases; see Shenfeld and van Handel \cite{SH18}, and \cite{HXZ21,Liu22,FHX23}.

\begin{remark}
We would like to mention that,  thanks to the anisotropic Heintze-Karcher inequality (cf. \cite{HLMG09}, and \cite[Prop.~6.1]{ALWX21}),  an extension of \autoref{Gauss1} to the setting of anisotropic geometry is achievable when $\partial_2\varphi=0$: 

Suppose $\varphi: (0,\infty)\to (0,\infty)$ is $C^1$-smooth with $\varphi'>0$.
Let $K, L$ be two smooth, strictly convex bodies with positive support functions $h_K,h_L$ and Gauss curvatures $\cK_K,\cK_L$. If
\eq{
\frac{\cK_L}{\cK_K}=\varphi\left(\frac{h_K}{h_L}\right),
}
then $K$ is a rescaling of $L$. 

The details are left to the interested reader.
\end{remark}

\section{background}
By a convex body in $\bbR^{n+1}$, we mean a compact convex set with non-empty interior. A smooth, strictly convex body is a convex body with $C^{\infty}$-smooth boundary and positive Gauss curvature.
 
Suppose $M=M^n\subset \bbR^{n+1}$ is a closed, smooth, strictly hypersurface with the outer unit normal vector field $\nu$ and support function $h$ and Gauss curvature $\cK$. Let $K$ denote the convex body with boundary $M$. The Gauss map of $M$, $\nu: M\to \bbS^n$, takes $x\in M$ to a unique point on $\bbS^n$ whose outer unit normal is $\nu(x)$. Here, we consider both $h$ and $\cK$ as functions on the unit sphere:
\eq{
h(x)=\langle \nu^{-1}(x),x\rangle,\quad 
\frac{1}{\cK(x)}=\frac{\det (\bar{\nabla}^2h+\bar{g}h)}{\det(\bar{g})}\Big|_x, \quad\quad x\in \bbS^n.
}

Let $\mu$ denote the spherical Lebesgue measure of the unit sphere and
 \eq{
  dV:=\frac{1}{n+1}\frac{h}{\cK}d\mu
  }
  be the cone volume measure of $K$, in particular, $V(K)=\int dV$. Note that $\cK^{-1}d\mu$ is the push-forward of the surface area measure of $M$, $d\mu_{M}$, via its Gauss map $\nu$.

Let $\cH$ denote the mean curvature of $M$ (i.e. the sum of the principal curvatures). By the Heintze-Karcher (cf. \cite{Ros87}), we have
\eq{\label{HK in GMP}
n\int \frac{1}{\cH}d\mu_{M}\geq (n+1)V(K).
}
The equality holds if and only if $M$ and $\bbS^n$ are homothetic.

\section{Uniqueness}
For the definition of a $k$-convex function see \cite[Def.~3.1]{GMTZ10}.
\begin{lemma}\label{magic sigma k ineq} Let $u\in C^2(\bbS^n)$ be $k$-convex.
 Then
\eq{
k\int |Du|^2u\sigma_kd\mu\leq \int u^2(\sigma_1\sigma_k-(k+1)\sigma_{k+1})d\mu+k\frac{|\int u\sigma_k Dud\mu|^2}{\int u\sigma_kd\mu}.
}
\end{lemma}
\pf{
The lemma when $u$ is the support function of a closed, smooth, strictly convex hypersurface was proved in \cite{IM23}. The proof in the $k$-convex case is similar; use \cite[Thm.~4.1]{GMTZ10} instead of \cite[Thm.~7.6.8]{Sch14}.
}

\begin{lemma}\label{HK reformulation}
We have
\eq{
\int \frac{\bar{\Delta}h}{\cK} d\mu\geq 0.
}
\end{lemma}
\pf{
Note that $\bar{\Delta}h+nh$ is the sum of the principal radii of curvature of $M$. Therefore,
$
\left(\bar{\Delta}h+nh\right)\cH\circ \nu^{-1}\geq n^2
$
and 
\eq{
\int \left(\bar{\Delta}h+nh\right)\frac{1}{\cK} d\mu\geq n^2\int \frac{1}{\cH}d\mu_{M}.
}
By the Heintze-Karcher inequality,
\eq{
\int \left(\bar{\Delta}h+nh\right)\frac{1}{\cK} d\mu\geq n(n+1)V(K)= n\int \frac{h}{\cK}d\mu.
}

}

\pf{[Proof of \autoref{sigma k soliton}]
In view of \autoref{magic sigma k ineq}, and $\int u\sigma_kDu=0,$ 
\eq{
k\alpha_k\int |Du|^2u\sigma_kd\mu\leq \alpha_k\int u^2\sigma_1\sigma_k-(k+1)u^2\sigma_{k+1}d\mu.
}
Recall that
$
(k+1)\sigma_{k+1}=\sigma_{k+1}^{ij}(\bar{\nabla}^2_{i,j}u+\bar{g}_{ij}u).
$
Suppose $\{e_i\}_{i=1}^n$ is an orthonormal basis of $T_{x_0}\Sn$ consisting of eigenvectors of $\bar{\nabla}^2u+\bar{g}u$ such that $(\bar{\nabla}^2u+\bar{g}u)(e_i,e_j)=\lambda_i\delta_{ij}$.
Note that $\bar{\nabla}_i\sigma_{k+1}^{ij}=0$, $e_i(Du)=\lambda_ie_i$, $\lambda_i\geq0$ and $\partial_2\varphi\geq 0$. Hence, by summing over $k$ and integrating by parts we find
\eq{
\int (u^2+|\bar{\nabla} u|^2)(k\alpha_ku\sigma_k)d\mu
&\leq \int nu^2\varphi-(\varphi+u\partial_1 \varphi)|\bar{\nabla} u|^2d\mu\\
&\hp{=}+\int -\alpha_ku^3\sigma_{k+1}^{ij}\bar{g}_{ij}+ 2\alpha_ku\sigma_{k+1}^{ij}\partial_iu\partial_ju d\mu.
}
Now due to the well-known identities
\eq{
\sigma_{k+1}^{ii}=\sigma_k-\lambda_i\sigma_k^{ii},\quad
\sigma_{k+1}^{ii}\bar{g}_{ii}=(n-k)\sigma_k,
}
we obtain
\eq{
\int \left(\sum_kk\alpha_ku\sigma_k-\varphi+u\partial_1\varphi\right)|\bar{\nabla} u|^2d\mu\leq 0.
}
Hence, $u$ is constant.
}

\pf{[Proof of \autoref{Gauss1}] By \autoref{HK reformulation}, we have
\eq{
\int \varphi(h,|Dh|)\bar{\Delta}hd\mu\geq 0.
}
Let $\{e_i\}_{i=1}^n$ be an orthonormal basis of $T_{x_0}\bbS^n$ consisting of eigenvectors of $\tau=\bar{\nabla}^2h+\bar{g}h$ and $\tau(e_i,e_j)=\lambda_i\delta_{ij}$. Since $e_i (Dh)=\lambda_i e_i$ and $\lambda_i>0$, the claim follows from integration by parts.
}

\providecommand{\bysame}{\leavevmode\hbox to3em{\hrulefill}\thinspace}
\providecommand{\href}[2]{#2}

\vspace{10mm}
\textsc{Institut f\"{u}r Diskrete Mathematik und Geometrie,\\ Technische Universit\"{a}t Wien, Wiedner Hauptstra{\ss}e 8-10,\\ 1040 Wien, Austria,} \email{\href{mailto:mohammad.ivaki@tuwien.ac.at}{mohammad.ivaki@tuwien.ac.at}}
\end{document}